\documentclass{llncs}
\usepackage{amssymb}
\usepackage{graphicx}

\usepackage{amsmath,epsfig}

\usepackage{url}
\urldef{\mailsa}\path|| \urldef{\mailsb}\path||
\urldef{\mailsc}\path||

\newtheorem{Assumption}{Assumption}

\catcode`\@=11
\def\mequal{\mathrel{\mathpalette\@mvereq{\hbox{\sevenrm m}}}}
\def\@mvereq#1#2{\lower.5\p@\vbox{\baselineskip\z@skip\lineskip1.5\p@
    \ialign{$\m@th#1\hfil##\hfil$\crcr#2\crcr=\crcr}}}
\def\partr#1#2{/\kern-.08333em/_{#1,#2}^{\phantom{.}}}
\def\invpartr#1#2{/\kern-.08333em/_{#1,#2}^{-1}}
\def\hpartr#1#2{/\kern-.08333em/_{#1,#2}^{h}}
\def\Epartr#1#2{/\kern-.08333em/_{#1,#2}^{E}}

\def\newdot{{\kern.8pt\cdot\kern.8pt}}
\def\,{\relax\ifmmode\mskip\thinmuskip\else\thinspace\fi}
\def\{{\relax\ifmmode\lbrace\else $\lbrace$\fi}
\def\}{\relax\ifmmode\rbrace\else $\rbrace$\fi}

\font\sevenrm=cmr7

\newcommand\DD{\mathbb{D}}
\newcommand\EE{\mathbb{E}}

\newcommand\RR{\mathbb{R}}

\def\grad{\mathpal{grad}}

\def\grad{\mathop{\rm grad}\nolimits}
\def\di{\displaystyle}
\def\f{\frac}
\def\a{\alpha }
\def\b{\beta }

\def\d{\delta }
\def\e{\varepsilon }

\def\l{\lambda }
\def\n{\nabla }

\def\s{\sigma }

\def\supp{\mathop{\rm supp}\nolimits}

\def\grad{\mathop{\rm grad}\nolimits}

\def\inj{\mathop{\rm inj}\nolimits}

\def\arcsinh{\mathop{\rm arcsinh}\nolimits}
\def\lip{\mathop{\rm Lip}\nolimits}
\def\sign{\mathop{\rm sign}\nolimits}

\makeatletter
\newcommand{\Rmnum}[1]{\expandafter\@slowromancap\romannumeral#1@}

\makeatother

\begin{document}
\index{Lastname, FirstName}
\mainmatter  

\title{Medians and means in Riemannian geometry: existence, uniqueness and computation}


\titlerunning{Lecture Notes in Computer Science}

%
%
\author{Marc Arnaudon\inst{1}, Fr\'{e}d\'{e}ric
Barbaresco\inst{2} and Le Yang\inst{1}
}
%

\institute{Laboratoire de Math\'{e}matiques et Applications, CNRS :
UMR 6086,\\ \quad\,Universit\'{e} de Poitiers,
86962 Futuroscope Chasseneuil Cedex, France\\
\and Thales Air Systems, Surface Radar, Technical Directorate,\\
Advanced Developments Dept. F-91470 Limours, France
}

%
%

\toctitle{Lecture Notes in Computer Science}
\tocauthor{Authors' Instructions}
\maketitle

\begin{abstract}
This paper is a short summary of our recent work on the medians and
means of probability measures in Riemannian manifolds. Firstly, the
existence and uniqueness results of local medians are given. In
order to compute medians in practical cases, we propose a
subgradient algorithm and prove its convergence. After that,
Fr\'{e}chet medians are considered. We prove their statistical
consistency and give some quantitative estimations of their
robustness with the aid of upper curvature bounds. We also show
that, in compact Riemannian manifolds, the Fr\'{e}chet medians of
generic data points are always unique. Stochastic and deterministic
algorithms are proposed for computing Riemannian $p$-means. The rate
of convergence and error estimates of these algorithms are also
obtained. Finally, we apply the medians and the Riemannian geometry
of Toeplitz covariance matrices to radar target detection.
\end{abstract}

\section{Introduction}
It has been widely accepted that the history of median begins from
the following question raised by P. Fermat in 1629: given a triangle
in the plan, find a point such that the sum of its distances to the
three vertices of the triangle is minimum. It is well known that the
answer to this question is: if each angle of the triangle is smaller
than $2\pi/3$, then the minimum point is such that the three
segments joining it and the vertices of the triangle form three
angles equal to $2\pi/3$; in the opposite case, the minimum point is
the vertex whose angle is no less than $2\pi/3$. This point is
called the median or the Fermat point of the triangle.

The notion of median also appears in statistics since a long time
ago. In 1774, when P. S. Laplace tried to find an appropriate notion
of the middle point for a group of observation values, he introduced
``the middle of probability'', the point that minimizes the sum of
its absolute differences to data points, this is exactly the one
dimensional median used by us nowadays.

A sufficiently general notion of median in metric spaces was
proposed in 1948 by M. Fr\'{e}chet in his famous article
\cite{Frechet}, where he defined a $p$-mean of a random variable $X$
to be a point which minimizes the expectation of its distance at the
power $p$ to $X$. This flexible definition allows us to define
various typical values, among which there are two important cases:
$p=1$ and $p=2$, corresponding to the notions of median and mean,
respectively.

Apparently, the median and mean are two notions of centrality for
data points. As a result, one may wonder that which one is more
advantageous? Statistically speaking, the answer to this question
depends on the distribution involved. For example, the mean has
obvious advantage over the median when normal distributions are
used. On the contrary, as far as Cauchy distributions are concerned,
the empirical mean has the same accuracy as one single observation,
so that it would be better to use the median instead of the mean in
this situation. Perhaps the most significant advantage of the median
over the mean is that the former is robust but the latter is not,
that is to say, the median is much less sensitive to outliers than
the mean. Roughly speaking, in order to move the median of a group
of data points to arbitrarily far, at least a half of data points
should be moved. Oppositely, in order to move the mean of a group of
data points to arbitrarily far, it suffices to move one data point.
So that medians are in some sense more prudent than means, as argued
by M. Fr\'{e}chet. The robustness property makes the median an
important estimator in situations when there are lots of noise and
disturbing factors.

The first formal definition of means for probability measures on
Riemannian manifolds was made by H. Karcher in \cite{Karcher}. To
introduce Karcher's result concerning means, consider a Riemannian
manifold $M$ with Riemannian distance $d$ and
\[B(a,\rho)=\{x\in M: d(x,a)<\rho\}\]
is a geodesic ball in $M$ centered at $a$ with a finite radius
$\rho$. Let $\Delta$ be an upper bound of sectional curvatures in
$\bar{B}(a,\rho)$ and $\inj$ be the injectivity radius of
$\bar{B}(a,\rho)$. Under the following condition:
\begin{equation}\label{radius condition}
\rho<\min\big\{\,\frac{\pi}{4\sqrt{\Delta}}\,,\,\frac{\inj}{2}\,\big\},
\end{equation}
where if $\Delta\leq0$, then $\pi/(4\sqrt{\Delta})$ is interpreted
as $+\infty$, Karcher showed that, with the aid of estimations of
Jacobi fields, the local energy functional
\begin{equation}\label{energy fun}
F_{\mu}:\quad \bar{B}(a,\rho)\longrightarrow \mathbf{R}_+,\quad
x\longmapsto\int_M d^2(x,p)\mu(dp)
\end{equation}
is strictly convex, thus it has a unique minimizer $b(\mu)$, which
he called the Riemannian center of mass of the probability measure
$\mu$. Moreover, $b(\mu)$ is also the unique solution of the
following equation:
\begin{equation}\label{exponential barycenter}
\int_{M}\exp_{x}^{-1}p\,\,\mu(dp)=0_{x},\quad x\in \bar{B}(a,\rho).
\end{equation}
From then on, local means of probability measures on Riemannian
manifolds are also called Karcher means, meanwhile, global means are
often called Fr\'{e}chet means. A rather general result concerning
the uniqueness of local means was proved by W. S. Kendall in
\cite{Kendall1}. As a particular case of Kendall's result, the
condition
\begin{equation}\label{Afsari condition}
\rho<\frac{1}{2}\min\big\{\,\frac{\pi}{\sqrt{\Delta}}\,,\,\inj\,\big\}
\end{equation}
is sufficient to ensure the uniqueness of the Kacher means of $\mu$.

Some generalizations of Karcher mean are given by many authors. For
instance, M. Emery and G. Mokobodzki defined in \cite{Emery} the
exponential barycenters and convex barycenters for measures on
affine manifolds. They also showed that a point $x$ is a convex
barycenter of a probability $\mu$ if and only if there exists a
continuous martingale starting from $x$ with terminal law $\mu$. The
uniqueness of exponential barycenters are generalized by M. Arnaudon
and X. M. Li in \cite{Arnaudon3} to probability measures on convex
affine manifolds with semilocal convex geometry. Moreover, the
behavior of exponential barycenters when measures are pushed by
stochastic flows is also considered in \cite{Arnaudon3}. In order to
study harmonic maps between Riemannian manifolds with probabilistic
methods, J. Picard also gave a generalized notion of barycenters in
\cite{Picard}. As we noted before, Karcher means are only local
minimizers of the energy functional $f_{\mu}$ in (\ref{energy fun}),
but it is easily seen that $f_{\mu}$ can be defined not only on the
closed ball $\bar{B}(a,\rho)$ but also on the whole manifold $M$ as
long as the second moment of $\mu$ is finite. This leads to the
global minimizers of the second moment function of $\mu$, which is
just the original definition of means made by Fr\'{e}chet. Global
minimizers are more useful in statistics than local ones, so that it
is necessary to know whether or under which conditions the Karcher
mean of $\mu$ is in fact the Fr\'{e}chet mean. For the case when
$\mu$ is a discrete measure supported by finitely many points in the
closed upper hemisphere, S. R. Buss and J. P. Fillmore showed in
\cite{Buss} that if the support of $\mu$ is not totally contained in
the equator then $\mu$ has a unique Karcher mean which lies in the
open hemisphere and equals to the Fr\'{e}chet mean. Inspired by the
methods of Buss and Fillmore, B. Afsari showed in \cite{Afsari} that
if the upper curvature bound $\Delta$ and the injectivity radius
$\inj$ in (\ref{Afsari condition}) is replaced by the ones of the
larger ball $B(a,2\rho)$, then all the Fr\'{e}chet $p$-means of
$\mu$ lie inside $B(a,\rho)$. Particularly, the Karcher mean
coincides with the Fr\'{e}chet mean. The existence and uniqueness of
$p$-means in Finsler geometry are recently proved by M. Arnaudon and
F. Nielsen in \cite{Arnaudon6}. They also showed that Finslerian
$p$-means are limiting points of continuous gradient flows and
developed algorithms for computing $p$-means in Finsler geometry.

Medians of discrete sample points on the sphere are studied by
economists and operational research experts in the 1970s and 1980s,
but they used the name ``location problems on a sphere''. For data
points lying in a spherical disc of radius smaller than $\pi/4$,
Drezner and Wesolowsky showed in \cite{Drezner78} that the cost
function is unimodal in that disc and the Fr\'{e}chet median is
unique if the data points are not contained in a single great
circle. It is also shown by Z. Drezner in \cite{Drezner81} that if
all the sample points are contained in a great circle, then one of
the sample points will be a Fr\'{e}chet median. Perhaps the first
work about Fr\'{e}chet medians on Riemannian manifolds is the paper
\cite{Noda} by R. Noda and his coauthors. They proved the
uniqueness, characterizations and position estimations of
Fr\'{e}chet medians for discrete sample points lying in a
Cartan-Hadamard manifold. In order to do robust statistics for data
living in a Riemannian manifold P. T. Fletcher and his coauthors
defined in \cite{Fletcher3} the local medians for discrete sample
points and showed their existence and uniqueness.

In this paper, we present our results on medians and means of
probability measures in Riemannian manifolds. Above all, the
motivation of our work: radar target detection is introduced in
section \ref{motivation}. After that, in section \ref{rme} we define
local medians for probability measures in Riemannian manifolds and
consider the problems of uniqueness and approximation. Under the
assumption that the support of the probability measure is contained
in a convex ball, we give some results on the characterization, the
position estimation and the uniqueness of medians. Then we propose a
subgradient algorithm to estimate medians as well as giving its
convergence result without condition of the sign of curvatures. Our
algorithm improves the one proposed in \cite{Fletcher3} which is
shown to be convergent only if the manifold is nonnegatively curved.
Finally, the problem of error estimation and rate of convergence are
also considered.

The aim of section \ref{Frechet medians} is to give some basic
properties of Fr\'{e}chet medians of probability measures in
Riemannian maniolds. Firstly, we give the consistency result of
Fr\'{e}chet medians in proper metric spaces. Particularly, if a
probability measure has only one Fr\'{e}chet median, then any
sequence of empirical Fr\'{e}chet medians will converge almost
surely to it. After that, we study the robustness of Fr\'{e}chet
medians in Riemannian manifolds. It is well known that in Euclidean
spaces, if a group of data points has more than a half concentrated
in a bounded region, then its Fr\'{e}chet median cannot be drown
arbitrarily far when the other points move. A generalization and
refinement of this result for data points in Riemannian manifolds is
given in Theorem \ref{1 robust}. This theorem also generalizes a
result in \cite{Afsari} which states that if the probability measure
is supported in a strongly convex ball, then all its Fr\'{e}chet
medians lie in that ball. At the end of this section, the uniqueness
question of Fr\'{e}chet sample medians is considered in the context
of compact Riemannian manifolds. It is shown that, apart from
several events of probability zero, the Fr\'{e}chet sample medians
are unique if the sample vector has a density with respect to the
canonical Lebesgue measure of the product manifold. In other words,
the Fr\'{e}chet medians of generic data points are always unique.

Section \ref{algo} is devoted to presenting algorithms for computing
Fr\'{e}chet $p$-means in order to meet practical needs.
Theorem~\ref{1 2.P1} gives stochastic algorithms which converge
almost surely to $p$-means in manifolds, which are easier to
implement than gradient descent algorithm since computing the
gradient of the function to minimize is not needed. The idea is at
each step to go in the direction of a point of the support of $\mu$.
The point is chosen at random according to $\mu$ and the size of the
step is a well chosen function of the distance to the point, $p$ and
the number of the step. The speed of convergence is given by
Theorem~\ref{1 2.P2}, which says that the renormalized inhomogeneous
Markov chain of Theorem~\ref{1 2.P1} converges in law to an
inhomogeneous diffusion process. We give the explicit expression of
this process, as well as its local characteristic. After that, the
performance of the stochastic algorithms are illustrated by
simulations. Finally, we show that the $p$-mean of $\mu$ can also be
computed by the method of gradient descent. The questions concerning
the choice of stepsizes and error estimates of this deterministic
method are also considered. We note that, for the case when
$p=+\infty$, M. Arnaudon and F. Nielsen developed in
\cite{Arnaudon5} an efficient algorithm to compute the circum-center
$e_{\infty}$ of probability measures in Riemannian manifolds.

In section \ref{Toeplitz}, we consider the manifold of $n\times n$
Toeplitz covariance matrices parameterized by the reflection
coefficients which are derived from Levinson's recursion of
autoregressive models. The explicit expression of the
reparametrization and its inverse are obtained. With the Riemannian
metric given by the Hessian of a K\"{a}hler potential, we show that
the manifold is in fact a Cartan-Hadamard manifold with lower
sectional curvature bound $-4$. After that, we compute the geodesics
and use the subgradient algorithm introduced in section \ref{rme} to
find the median of Toeplitz covariance matrices. Finally, we give
some simulated examples to illustrate the application of the median
method to radar target detection.

\section{Motivation: radar target detection}\label{motivation}

Suggested by J. C. Maxwell's seminal work on electromagnetism, H.
Hertz carried out an experiment in 1886 which validated that radio
waves could be reflected by metallic objects. This provided C.
H\"{u}elsmeyer the theoretical foundation of his famous patent on
``telemobiloscope'' in 1904. He showed publicly in Germany and
Netherlands that his device was able to detect remote metallic
objects such as ships, even in dense fog or darkness, so that
collisions could be avoided. H\"{u}elsmeyer's ``telemobiloscope'' is
recognized as the primogenitor of modern radar even though it could
only detect the direction of an object, neither its distance nor its
speed. This is because the basic idea of radar was already born:
send radio waves in a predetermined direction and then receive the
possible echoes reflected by a target. In order to know the distance
and the radial speed of the target, it suffices to send successively
two radio waves. In fact, it is easily seen that the distance $d$ of
the target can be computed by the formula
$$d=\frac{c\Delta t}{2},$$
where $c$ is the speed of light and $\Delta t$ is the time interval
between every emission and reception in the direction under test.
Moreover, the radial speed $v$ of the target can be deduced by the
Doppler effect which states that the frequency of a wave is changed
for an observer moving relatively to the source of the wave. More
precisely,
$$v=\frac{\lambda\Delta\varphi}{4\pi\Delta t},$$
where $\lambda$ and $\Delta\varphi$ are the wavelength and the
skewing of the two emitted radio waves, respectively. As a result,
the direction, the distance and the speed of the target can all be
determined.

For simplicity, from now on we only consider a fixed direction in
which a radar sends radio waves. Since the range of emitted waves
are finite, we can divide this direction into some intervals each of
which represents a radar cell under test. The radar sends each time
a rafale of radio waves in this direction and then receive the
returning echoes. For each echo we measure its amplitude $r$ and
phase $\varphi$, so that it can be represented by a complex number
$z=re^{i\varphi}$. As a result, the observation value of each radar
cell is a complex vector $Z=(z_1,\dots,z_N)$, where $N$ is the
number of waves emitted in each rafale.

The aim of target detection is to know whether there is a target at
the location of some radar cell in this direction. Intuitively
speaking, a target is an object whose behavior on reflectivity or on
speed is very different from its environment. The classical methods
for target detection is to compute the difference between the
discrete Fourier transforms of the radar observation values of the
cell under test and that of its ambient cells. The bigger this
difference is, the more likely a target appears at the location of
the cell under test. However, the performance of these classical
methods based on Doppler filtering using discrete Fourier transforms
together with the Constant False Alarm Rate (CFAR) is not very
satisfactory due to their low resolutions issues in perturbed radar
environment or with smaller bunch of pulses.

In order to overcome these drawbacks, a lot of mathematical models
for spectra estimation were introduced, among which the method based
on autoregressive models proposed by F. Barbaresco in
\cite{Barbaresco1} is proved to be very preferable. We shall
introduce this method in Chapter \ref{Toeplitz} of this
dissertation. The main difference between this new method and the
classical ones is that, instead of using directly the radar
observation value $Z$ of each cell, we regard it as a realization of
a centered stationary Gaussian process and identify it to its
covariance matrix $R=\mathbf{E}[ZZ^*]$. Thus the new observation
value for each radar cell is a covariance matrix which is also
Toeplitz due to the stationarity of the process. As a result, the
principle for target detection becomes to find the cells where the
covariance matrix differs greatly from the average matrix of its
neighborhood. Once such cells are determined we can conclude that
there are targets in these locations. In order to carry out this new
method, there are two important things which should be considered
seriously. One is to define a good distance between two Toeplitz
covariance matrices. The other is to give a reasonable definition of
the average of covariance matrices, which should be robust to
outliers so as to be adapted to perturbed radar environment, and
develop an efficient method to compute it in practical cases. These
works will be done in the following by studying the Riemannian
geometry of Toeplitz covariance matrices and the medians of
probability measures in Riemannian manifolds.

\section{Riemannian median and its estimation}\label{rme}

In this section, we define local medians of a probability measure on
a Riemannian manifold, give their characterization and a natural
condition to ensure their uniqueness. In order to compute medians in
practical cases, we also propose a subgradient algorithm and show
its convergence. The mathematical details of this section can be
found in \cite{Yang}.

In more detail, let $M$ be a complete Riemannian manifold with
Riemannian metric $\langle\,\cdot\,,\cdot\,\rangle$ and Riemannian
distance $d$. We fix an open geodesic ball
\[B(a,\rho)=\{x\in M: d(x,a)<\rho\}\] in $M$ centered at $a$ with a
finite radius $\rho$. Let $\delta$ and $\Delta$ denote respectively
a lower and an upper bound of sectional curvatures $K$ in
$\bar{B}(a,\rho)$. The injectivity radius of $\bar{B}(a,\rho)$ is
denoted by $\inj\,(\bar{B}(a,\rho))$. Furthermore, we assume that
the radius of the ball verifies
\begin{equation}\label{radius condition}
\rho<\min\big\{\,\frac{\pi}{4\sqrt{\Delta}}\,,\,\frac{\inj\,(\bar{B}(a,\rho))}{2}\,\big\},
\end{equation}
where if $\Delta\leq0$, then $\pi/(4\sqrt{\Delta})$ is interpreted
as $+\infty$.

We consider a probability  measure $\mu$ on $M$ whose support is
contained in the open ball $B(a,\rho)$ and define a function
\[f:\qquad \bar{B}(a,\rho) \longrightarrow\mathbf{R}_+\,,\qquad x\longmapsto\int_M d(x,p)\mu(dp).\]

This function is 1-Lipschitz, hence continuous on the compact set
$\bar{B}(a,\rho)$. The convexity of the distance function on
$\bar{B}(a,\rho)$ yields that $f$ is also convex. Hence we don't
need to distinguish its local minima from its global ones. Now we
can give the following definition:

\begin{definition}
A minimum point of $f$ is called a median of $\mu$. The set of all
the medians of $\mu$ will be denoted by $\mathfrak{M}_{\mu}$. The
minimal value of $f$ will be denoted by $f_*$.
\end{definition}

It is easily seen that $\mathfrak{M}_{\mu}$ is compact and convex.
Moreover, by computing the right derivative of $f$ we can prove the
following characterization of $\mathfrak{M}_{\mu}$.

\begin{theorem}\label{1 characterization of median}
The set $\mathfrak{M}_{\mu}$ is characterized by
$$\mathfrak{M}_{\mu}=\big\{x\in
\bar{B}(a,\rho):|H(x)|\leq\mu{\{x\}}\big\},$$ where for $x\in
\bar{B}(a,\rho)$,
\[H(x):=\int_{M\setminus\{x\}}\frac{-\exp_x^{-1}p}{d(x,p)}\mu(dp),\]
is a tangent vector at $x$ satisfying $|H(x)|\leq 1$.
\end{theorem}

Observing that every geodesic triangle in $\bar{B}(a,\rho)$ has at
most one obtuse angle, we can prove the following result which gives
a position estimation for the medians of $\mu$.

\begin{proposition}\label{1 location}
$\mathfrak{M}_{\mu}$ is contained in the smallest closed convex
subset of $B(a,\rho)$ containing the support of $\mu$.
\end{proposition}

In Euclidean case, it is well known that if the sample points are
not collinear, then their medians are unique. Hence we get a natural
condition of $\mu$ to ensure the uniqueness for medians in Riemannian case:\\

$\ast$ \quad \emph{The support of $\mu$ is not totally contained in
any geodesic. This means that for every geodesic $\gamma$:
$[\,0,1\,]\rightarrow \bar{B}(a,\rho)$, we have
$\mu(\gamma[\,0,1\,])<1$.}\\

This condition implies that $f$ is strictly convex along every
geodesic in $\bar{B}(a,\rho)$, so that it has one and only one
minimizer, as stated by the theorem below.

\begin{theorem}\label{1 uniqueness of meidan}
If condition $\ast$ holds, then $\mu$ has a unique median.
\end{theorem}

With further analysis, we can show a stronger quantitative version
of Theorem \ref{1 uniqueness of meidan}, which is crucial in the
error estimations of the subgradient algorithm as well as in the
convergence proof of the stochastic algorithm for computing medians
in section \ref{algo}.

\begin{theorem}\label{1 strong convex}
If condition $\ast$ holds, then there exits a constant $\tau>0$ such
that for every $x\in\bar{B}(a,\rho)$ one has
\[f(x)\geq f_{\ast}+\tau d^2(x,m),\]
where $m$ is the unique median of $\mu$.
\end{theorem}

The main results of approximating medians of $\mu$ by subgradient
method is summarized in the following theorem. The idea stems from
the basic observation that $H(x)$ is a subgradient of $f$ at $x$ for
every $x\in\bar{B}(a,\rho)$.

\begin{theorem}\label{1 convergence}
Let $(t_k)_k$ be a sequence of real numbers such that
\[t_k>0,\quad\lim_{k\rightarrow\infty}t_k=0\quad\text{and}\quad\sum_{k=0}^{\infty}t_k=+\infty.\]
Define a sequence $(x_k)_k$ by $x_0\in\bar{B}(a,\rho)$ and for
$k\geq0$,
\[
x_{k+1}=
\begin{cases}
x_k, & \text{if $H(x_k)=0$};\\
\exp_{x_k}\bigg(-t_k\cfrac{H(x_k)}{|H(x_k)|}~\bigg),  & \text{if
$H(x_k)\neq0$.}
\end{cases}
\]
Then there exists some constant $T>0$ such that if we choose
$t_k\leq T$ for every $k\geq0$, then the sequence $(x_k)_k$ is
contained in $\bar{B}(a,\rho)$ and verifies
\[\lim_{k\rightarrow\infty}d(x_k,\mathfrak{M}_{\mu})=0\quad\text{and}\quad\lim_{k\rightarrow\infty}f(x_k)=f_*.\]
Moreover, if the sequence $(t_k)_k$ also verifies
\[\sum_{k=0}^{\infty}t^2_k<+\infty,\]
then there exists some $m\in \mathfrak{M}_{\mu}$ such that
$x_k\longrightarrow m$.
\end{theorem}

\begin{remark}
We can choose the constant $T$ in Theorem \ref{1 convergence} to be
\[T=\frac{\rho-\sigma}{C(\rho,\delta)F(\rho,\Delta)+1},\]
where $\sigma=\sup\{d(p,a):p\in \supp\mu\}$,
\[
F(\rho,\Delta)=
\begin{cases}
1,  & \text{if $\Delta\geq0$;}\\

\cosh(2\rho\sqrt{-\Delta}), & \text{if $\Delta<0$,}
\end{cases}
\]
and
\[
C(\rho,\delta)=
\begin{cases}
1,  & \text{if $\delta\geq0$;}\\

2\rho\sqrt{-\delta}\coth(2\rho\sqrt{-\delta}), & \text{if
$\delta<0$.}
\end{cases}
\]
\end{remark}

The proposition below gives the error estimation of the algorithm in
Theorem \ref{1 convergence}.

\begin{proposition}\label{1 error}
Let condition $\ast$ hold and the stepsizes $(t_k)_k$ in Theorem
\ref{1 convergence} satisfy
 \[\lim_{k\rightarrow\infty}t_k=0\quad \text{and}\quad\sum_{k=0}^{\infty}t_k=+\infty.\]
Then there exists $N\in\mathbf{N}$, such that for every $k\geq N$,
\[d^2(x_k,m)\leq b_k,\]
where $m$ is the unique median of $\mu$ and the sequence
$(b_k)_{k\geq N}$ is defined by
\[b_N=(\rho+\sigma)^2\quad\text{and}\quad b_{k+1}=(1-2\tau
t_k)b_k+C(\rho,\delta)t^2_k\,,\quad k\geq N,\] which converges to
$0$ when $k\rightarrow\infty$. More explicitly, for every $k\geq N$,
\[b_{k+1}=(\rho+\sigma)^2\prod_{i=N}^{k}(1-2\tau t_i)
+C(\rho,\delta)\big(\sum_{j=N+1}^{k}t^2_{j-1}\prod_{i=j}^{k}(1-2\tau
t_i)+t^2_k\,\,\big).\]
\end{proposition}

\section{Some properties of Fr\'{e}chet medians in Riemannian
manifolds}\label{Frechet medians}

This section is devoted to some basic results about Fr\'{e}chet
medians, or equivalently, global medians. We show the consistency of
Fr\'{e}chet medians in proper metric spaces, give a quantitative
estimation for the robustness of Fr\'{e}chet medians in Riemannian
manifolds and show the almost sure uniqueness of Fr\'{e}chet sample
medians in compact Riemannian manifolds. We refer to \cite{Yang3}
for more details of this section.

\subsection{Consistency of Fr\'{e}chet medians in metric spaces}

In this subsection, we work in a proper metric space $(M,d)$ (recall
that a metric space is proper if and only if every bounded and
closed subset is compact). Let $P_1(M)$ denote the set of all the
probability measures $\mu$ on $M$ verifying
 \[\int_M d(x_0,p)\mu(dp)<\infty,\,\,\text{for\,\,some}\,\,x_0\in M.\]
For every $\mu\in P_1(M)$ we can define a function
\[f_{\mu}:\qquad M \longrightarrow\mathbf{R}_+\,,\qquad x\longmapsto\int_M d(x,p)\mu(dp).\]
This function is 1-Lipschitz hence continuous on $M$. Since $M$ is
proper, $f_{\mu}$ attains its minimum (see \cite[p. 42]{Sahib}), 
so we can give the following definition:

\begin{definition}
Let $\mu$ be a probability measure in $P_1(M)$, then a \emph{global}
minimum point of $f_{\mu}$ is called a Fr\'{e}chet median of $\mu$.
The set of all the Fr\'{e}chet medians of $\mu$ is denoted by
$Q_{\mu}$. Let $f^{*}_{\mu}$ denote the global minimum of $f_{\mu}$.
\end{definition}

By the Kantorovich-Rubinstein duality of $L^1$-Wasserstein distance
(see \cite[p. 107]{Villani}), we can show that Fr\'{e}chet medians
are characterized by $1$-Lipschitz functions. A corresponding result
that Riemannian barycenters are characterized by convex functions
can be found in \cite[Lemma 7.2]{Kendall1}.

\begin{proposition}\label{1 lip}
Let $\mu\in P_1(M)$ and $M$ be also separable, then
$$Q_{\mu}=\bigg\{x\in M: \varphi(x)\leq f^*_{\mu}+\int_M \varphi(p)\mu(dp),\,\,\text{for every}\,\,\varphi\in
\lip_1(M)\bigg\},$$ where $\lip_1(M)$ denotes the set of all the
1-Lipschitz functions on $M$.

\end{proposition}

The following theorem states that the uniform convergence of first
moment functions yields the convergence of Fr\'{e}chet medians.

\begin{theorem}\label{1 consistency}
Let $(\mu_n)_{n\in\mathbf{N}}$ be a sequence in $P_1(M)$ and $\mu$
be another probability measure in $P_1(M)$. If $(f_{\mu_n})_n$
converges uniformly on $M$ to $f_{\mu}$, then for every
$\varepsilon>0$, there exists $N\in\mathbf{N}$, such that for every
$n\geq N$ we have
\[Q_{\mu_n}\subset B(Q_{\mu},\varepsilon):=\{x\in M: d(x,
Q_{\mu})<\varepsilon\}.\]
\end{theorem}

As a corollary to Theorem \ref{1 consistency}, Fr\'{e}chet medians
are strongly consistent estimators. The consistency of Fr\'{e}chet
means is proved in \cite{Bhattacharya}.

\begin{corollary}
Let $(X_n)_{n\in\mathbf{N}}$ be a sequence of i.i.d random variables
of law $\mu\in P_1(M)$ and $(m_n)_{n\in\mathbf{N}}$ be a sequence of
random variables such that $m_n\in Q_{\mu_n}$ with
$\mu_n=\frac{1}{n}\sum_{k=1}^{n}\delta_{X_k}$. If $\mu$ has a unique
Fr\'{e}chet median $m$, then $m_n\longrightarrow m$ a.s.
\end{corollary}

\subsection{Robustness of Fr\'{e}chet medians in Riemannian manifolds}

The framework of this subsection is a complete Riemannian manifold
$(M,d)$ whose dimension is no less than 2. We fix a closed geodesic
ball
\[\bar{B}(a,\rho)=\{x\in M: d(x,a)\leq\rho\}\] in $M$ centered at
$a$ with a finite radius $\rho>0$ and a probability measure $\mu\in
P_1(M)$ such that
\[\mu (\bar{B}(a,\rho))\triangleq\alpha>\frac{1}{2}.\]

The aim of this subsection is to estimate the positions of the
Fr\'{e}chet medians of $\mu$, which gives a quantitative estimation
for robustness. To this end, the following type of functions are of
fundamental importance for our methods. Let $x,z\in M$, define
\[h_{x,z}:\quad
\bar{B}(a,\rho)\longrightarrow\mathbf{R},\quad p\longmapsto
d(x,p)-d(z,p).\] Obviously, $h_{x,z}$ is continuous and attains its
minimum.

By a simple estimation on the minimum of $h_{x,a}$ we get the
following basic result.

\begin{theorem}\label{1 general estimation}
The set $Q_{\mu}$ of all the Fr\'{e}chet medians of $\mu$ verifies
\[Q_{\mu}\subset \bar{B}\bigg(a,\frac{2\alpha\rho}{2\alpha-1}\bigg):= B_*.\]
\end{theorem}

\begin{remark}
It is easily seen that the conclusion of Theorem \ref{1 general
estimation} also holds if $M$ is only a proper metric space.
\end{remark}

\begin{remark}
As a direct corollary to Theorem \ref{1 general estimation}, if
$\mu$ is a probability measure in $P_1(M)$ such that for some point
$m\in M$ one has $\mu\{m\}>1/2$, then $m$ is the unique Fr\'{e}chet
median of $\mu$.
\end{remark}

In view of Theorem \ref{1 general estimation}, let $\Delta$ be an
upper bound of sectional curvatures in $B_*$ and $\inj$ be the
injectivity radius of $B_*$. By computing the minima of some typical
functions $h_{x,z}$ in model spaces $\mathbb{S}^2$, $\mathbb{E}^2$
and $\mathbb{H}^2$,  and then comparing with the ones in $M$, we get
the following main result of this subsection.

\begin{theorem}\label{1 robust}
Assume that
\begin{equation}\label{1 assumption}
\frac{2\alpha\rho}{2\alpha-1}<r_*:=\min\{\frac{\pi}{\sqrt{\Delta}}\,\,,\inj\,\},
\end{equation}
where if $\Delta\leq0$, then $\pi/\sqrt{\Delta}$ is interpreted as
$+\infty$.\\

i) If $\Delta>0$ and $Q_{\mu}\subset \bar{B}(a,r_*/2)$, then
$$
Q_{\mu}\subset\bar{B}\bigg(a,\frac{1}{\sqrt{\Delta}}\arcsin\big(\frac{\alpha\sin(\sqrt{\Delta}\rho)}{\sqrt{2\alpha-1}}\big)\bigg).
$$
Moreover, any of the two conditions below implies $Q_{\mu}\subset
\bar{B}(a,r_*/2)$:

$$a)\quad\frac{2\alpha\rho}{2\alpha-1}\leq\frac{r_*}{2};
\qquad b)\quad
\frac{2\alpha\rho}{2\alpha-1}>\frac{r_*}{2}\quad\text{and}\quad
F_{\alpha,\rho,\Delta}(\frac{r_*}{2}-\rho)\leq0,$$ where $
F_{\alpha,\rho,\Delta}(t)=\cot(\sqrt{\Delta}(2\alpha-1)t)-\cot(\sqrt{\Delta}
t)-2\cot(\sqrt{\Delta}\rho),\,\, t\in(0,\cfrac{\rho}{2\alpha-1}\,\,].$\\

ii) If $\Delta=0$, then
$$Q_{\mu}\subset\bar{B}\bigg(a,\frac{\alpha\rho}{\sqrt{2\alpha-1}}\bigg).
$$

iii) If $\Delta<0$, then
$$Q_{\mu}\subset\bar{B}\bigg(a,\frac{1}{\sqrt{-\Delta}}\arcsinh\big(\frac{\alpha\sinh(\sqrt{-\Delta}\rho)}{\sqrt{2\alpha-1}}\big)\bigg).
$$
Finally any of the above three closed balls is contained in the open
ball $B(a,r_*/2)$.
\end{theorem}

\begin{remark}
Although we have chosen the framework of this section to be a
Riemannian manifold, the essential tool that has been used is the
hinge version of the triangle comparison theorem. Consequently,
Theorem \ref{1 robust}  remains true if $M$ is a CAT$(\Delta)$ space
(see \cite[Chapter 2]{Bridson}) and $r_*$ is replaced by
$\pi/\sqrt{\Delta}$.
\end{remark}

\begin{remark}
For the case when $\alpha=1$, the assumption \eqref{1 assumption}
becomes
$$\rho<\frac{1}{2}\min\{\frac{\pi}{\sqrt{\Delta}}\,\,,\inj\,\}.$$
Observe that in this case, when $\Delta>0$, the condition
$F_{1,\rho,\Delta}(r_*/2-\rho)\leq0$ is trivially true in case of
need. Hence Theorem \ref{1 robust} yields that
$Q_{\mu}\subset\bar{B}(a,\rho)$, which is exactly what the Theorem
2.1 in \cite{Afsari} says for medians.
\end{remark}

\subsection{Uniqueness of Fr\'{e}chet sample medians in
compact Riemannian manifolds}

Before introducing the results of this subsection we give some
notations. For each point $x\in M$, $S_x$ denotes the unit sphere in
$T_xM$. Moreover, for a tangent vector $v\in S_x$, the distance
between $x$ and its cut point along the geodesic starting from $x$
with velocity $v$ is denoted by $\tau(v)$. Certainly, if there is no
cut point along this geodesic, then we define $\tau(v)=+\infty$. For
every point $(x_1,\dots,x_N)\in M^N$, where $N\geq3$ is a fixed
natural number, we write
$$\mu(x_1,\dots,x_N)=\frac{1}{N}\sum_{k=1}^{N}\delta_{x_k}.$$
The set of all the Fr\'{e}chet medians of $\mu(x_1,\dots,x_N)$, is
denoted by $Q(x_1,\dots,x_N)$.

The following theorem states that in order to get the uniqueness of
Fr\'{e}chet medians, it suffices to move two data points towards a
common median along some minimizing geodesics for a little distance.

\begin{theorem}
Let $(x_1,\dots,x_N)\in M^N$ and $m\in Q(x_1,\dots,x_N)$. Fix two
normal geodesics $\gamma_1,\gamma_2:[0,+\infty)\rightarrow M$ such
that $\gamma_1(0)=x_1$, $\gamma_1(d(x_1,m))=m$, $\gamma_2(0)=x_2$
and $\gamma_2(d(x_2,m))=m$. Assume that
\[
x_2\notin
\begin{cases}
\gamma_1[0,\tau(\dot{\gamma}_1(0))],& \text{if $\tau(\dot{\gamma}_1(0))<+\infty$;}\\
\gamma_1[0,+\infty),& \text{if $\tau(\dot{\gamma}_1(0))=+\infty$.}
\end{cases}
\]
Then for every $t\in(0,d(x_1,m)]$ and $s\in(0,d(x_2,m)]$ we have
$$Q(\gamma_1(t),\gamma_2(s),x_3,\dots,x_N)=\{m\}.$$
\end{theorem}

Generally speaking, the non uniqueness of Fr\'{e}chet medians is due
to some symmetric properties of data points. As a result, generic
data points should have a unique Fr\'{e}chet median. In mathematical
language, this means that the set of all the particular positions of
data points is of Lebesgue measure zero. After eliminate all these
particular cases we obtain the following main result:

\begin{theorem}\label{1 unique}
Assume that $M$ is compact. Then $\mu(x_1,\dots,x_N)$ has a unique
Fr\'{e}chet median for almost every $(x_1,\dots,x_N)\in M^N$.
\end{theorem}

\begin{remark}
In probability language, Theorem \ref{1 unique} is equivalent to say
that if $(X_1,\dots,$\\$X_N)$ is an $M^N$-valued random variable
with density, then $\mu(X_1,\dots,X_N)$ has a unique Fr\'{e}chet
median almost surely. Clearly, the same statement is also true if
$X_1,\dots,X_N$ are independent and $M$-valued random variables with
density.
\end{remark}

\section{Stochastic and deterministic algorithms for computing means
of probability measures}\label{algo}

In this section, we consider a probability measure $\mu$ supported
by a regular geodesic ball in a manifold and, for any $p\ge 1$,
define a stochastic algorithm which converges almost surely to the
$p$-mean $e_p$ of $\mu$. Assuming furthermore that the functional to
minimize is regular around $e_p$, we prove that a natural
renormalization of the inhomogeneous Markov chain converges in law
into an inhomogeneous diffusion process. We give the explicit
expression of this process, as well as its local characteristic.
After that, the performance of the stochastic algorithms are
illustrated by simulations. Finally, we show that the $p$-mean of
$\mu$ can also be computed by the method of gradient descent. The
questions concerning the choice of stepsizes and error estimates of
this deterministic method are also considered. For more mathematical
details of this section, see \cite{Arnaudon4} and \cite{Yang4}.

\subsection{Stochastic algorithms for computing $p$-means}
Let $M$ be a Riemannian manifold whose sectional curvatures
$K(\sigma)$ verify $-\b^2\leq K(\sigma)\leq\a^2$, where $\a,\b$ are
positive numbers. Denote by $\rho$ the Riemannian distance on $M$.
Let $B(a,r)$ be a geodesic ball in $M$ and $\mu$ be a probability
measure with support included in a compact convex subset $K_\mu$ of
$B(a,r)$. Fix $p\in[1,\infty)$. We will always make the following
assumptions on $(r,p,\mu)$:
\begin{Assumption}
\label{1 A1} The support of $\mu$ is not reduced to one point.
Either $p>1$ or the support of $\mu$ is
 not contained in a line, and the radius $r$ satisfies
\begin{equation}\notag
 \label{rp}
r<r_{\a,p}\quad \hbox{with}\left\{\begin{array}{ccc}
                              r_{\a,p}&=\f12\min\left\{{\rm inj}(M),\f{\pi}{2\a}\right\}, &\hbox{if}\  p\in [1,2);\\
r_{\a,p}&=\f12\min\left\{{\rm inj}(M),\f{\pi}{\a}\right\},
&\hbox{if} \ p\in [2,\infty).
                             \end{array}
\right.
\end{equation}
 \end{Assumption}
Under Assumption~\ref{1 A1}, it has been proved in
\cite[Theorem~2.1]{Afsari}  that the function
\begin{equation}\notag
\label{2.1}
\begin{split}
H_p : M&\longrightarrow \RR_+\\
x&\longmapsto \int_M\rho^p(x,y)\mu(dy)
\end{split}
\end{equation}
has a unique minimizer $e_p$ in $M$, the $p$-mean of $\mu$, and
moreover $e_p\in B(a,r)$.  If $p=1$, $e_1$ is the median of $\mu$.

\begin{remark}
The existence and uniqueness of $p$-means in Finsler geometry are
recently proved by M. Arnaudon and F. Nielsen in \cite{Arnaudon6}.
They also showed that Finslerian $p$-means are limiting points of
continuous gradient flows and developed algorithms for computing
$p$-means in Finsler geometry.
\end{remark}

In the following theorem, we define a stochastic gradient algorithm
$(X_k)_{k\geq 0}$ to approximate the $p$-mean $e_p$ and prove its
convergence. In the sequel, let
\begin{equation}\notag
\label{K} K=\bar B(a, r-\e)\quad \hbox{with}\quad \e=
\f{\rho(K_\mu,B(a,r)^c)}{2}.
\end{equation}

\begin{theorem}
\label{1 2.P1} Let $(P_k)_{k\ge 1}$ be a sequence of independent
$B(a,r)$-valued random variables, with law $\mu$. Let $(t_k)_{k\ge
1}$ be a sequence of positive numbers satisfying
\begin{equation}\notag
\label{2.2bis} \forall k\ge 1,\ \ t_k\leq\min\left(\f1{C_{p,\mu,K}},
\f{\rho(K_\mu,B(a,r)^c)}{2p(2r)^{p-1}}\right),
\end{equation}
\begin{equation}\notag
\label{2.2} \sum_{k=1}^\infty t_k=+\infty\quad\hbox{and}\quad
\sum_{k=1}^\infty t_k^2<\infty,
\end{equation}
where $C_{p,\mu,K}>0$ is a constant.

Letting $x_0\in K$, define inductively the random walk $(X_k)_{k\ge
0}$ by
\begin{equation}\notag
\label{2.3} X_0=x_0\quad\hbox{and for $k\ge 0$}\quad
X_{k+1}=\exp_{X_{k}}\left(-t_{k+1}\grad_{X_{k}} F_p(\cdot,
P_{k+1})\right)
\end{equation}
where $F_p(x,y)=\rho^p(x,y)$, with the convention $\grad_x
F_p(\cdot, x)=0$.

The random walk $(X_k)_{k\ge 1}$ converges in $L^2$ and almost
surely to $e_p$.
\end{theorem}

\begin{remark}
For the case when $p=+\infty$, M. Arnaudon and F. Nielsen developed
in \cite{Arnaudon5} an efficient algorithm to compute the
circum-center $e_{\infty}$ of probability measures in Riemannian
manifolds.
\end{remark}

In the following example, we focus on the case $M=\RR^d$ and $p=2$
where drastic simplifications occur.
\begin{example}
\label{E1} In the case when $M=\RR^d$ and $\mu$ is a compactly
supported probability measure on $\RR^d$, the stochastic gradient
algorithm \eqref{2.3} simplifies into
\[
X_0=x_0\quad\hbox{and for $k\ge 0$}\quad
X_{k+1}=X_k-t_{k+1}\grad_{X_{k}} F_p(\cdot, P_{k+1}).
\]
If furthermore $p=2$, clearly $e_2=\EE[P_1]$ and $\grad_x F_p(\cdot,
y)=2(x-y)$, so that the linear relation
\[
X_{k+1}=(1-2t_{k+1})X_k+2t_{k+1}P_{k+1},\quad k\geq 0
\]
holds true and an easy induction proves that
\begin{equation} \label{2.12}
X_k=x_0\prod_{j=0}^{k-1}(1-2t_{k-j})+2\sum_{j=0}^{k-1}P_{k-j}t_{k-j}\prod_{\ell=0}^{j-1}(1-2t_{k-\ell}),\quad
k\ge 1.
\end{equation}
Now, taking $\di t_k=\f1{2k}$, we have
$$
\prod_{j=0}^{k-1}(1-2t_{k-j})=0\quad \hbox{and} \quad
\prod_{\ell=0}^{j-1}(1-2t_{k-\ell})=\f{k-j}{k}
$$
so that
$$
X_k=\sum_{j=0}^{k-1}P_{k-j}\f1k=\f1k\sum_{j=1}^kP_j.
$$
The stochastic gradient algorithm estimating the mean $e_2$ of $\mu$
is given by the empirical mean of a growing sample of independent
random variables with distribution $\mu$. In this simple case, the
result of Theorem~\ref{1 2.P1} is nothing but the strong law of
large numbers. Moreover, fluctuations around the mean are given by
the central limit theorem and Donsker's theorem.
\end{example}

The fluctuation of the random walk $(X_k)_k$ defined in Theorem
\ref{1 2.P1} is summarized in the following theorem.

\begin{theorem}\label{1 2.P2}
Assume that in Theorem \ref{1 2.P1}
$$
t_k=\min\left(\frac{\delta}{k}, \min\left(\f1{C_{p,\mu,K}},
\f{\rho(K_\mu,B(a,r)^c)}{2p(2r)^{p-1}}\right)\right), \quad k\geq 1,
$$
for some $\delta>0$. We define for $n\ge 1$ the Markov chain
$(Y_k^n)_{k\ge 0}$ in $T_{e_p}M$ by
$$
 \label{3.4}
Y_k^n=\f{k}{\sqrt{n}}\exp_{e_p}^{-1}X_k.
$$
Assume that $H_p$ is $C^2$ in a neighborhood of ~$e_p$ and
$\d>C^{-1}_{p,\mu,K}$. Then the sequence of processes $\di
\left(Y_{[nt]}^n\right)_{t\ge 0}$ converges weakly in
$\DD((0,\infty),T_{e_p}M)$ to a diffusion process given by

$$
\label{3.32}
 y_\d(t)=\sum_{i=1}^dt^{1-\d\l_i}\int_0^ts^{\d\l_i-1}\langle \d\s \,dB_s,e_i\rangle e_i,\quad t\ge
 0,
$$
where $B_t$ is the standard Brownian motion in $T_{e_p}M$ and $\s\in
{\rm End}(T_{e_p}M)$ satisfying
\[\s\s^\ast=\EE\left[\grad_{e_p}F_p(\cdot,P_1)\otimes
\grad_{e_p}F_p(\cdot, P_1)\right],
\]
$(e_i)_{1\le i\le d}$ is an orthonormal basis diagonalizing the
symmetric bilinear form $\n dH_p(e_p)$ and $(\l_i)_{1\le i\le d}$
are the associated eigenvalues.
\end{theorem}

\subsection{Simulations of stochastic algorithms}
\subsubsection{A non uniform measure on the unit square in the plane}
Here $M$ is the Euclidean plane $\RR^2$ and $\mu$ is the
renormalized restriction to the square $[0,4]\times [0,4]$ of an
exponential law on $[0,\infty)\times [0,\infty)$. The red path
represents one trajectory of the inhomogeneous Markov chain
$(X_k)_{k\ge 0}$ corresponding to $p=1$, with linear interpolation
between the different steps. The red point is $e_1$. Black circles
represent the values of $(P_k)_{k\ge 1}$.
\begin{figure}[!htpb]
\begin{center}
\includegraphics[bb=-8 -9 116 116, width =2.5in]{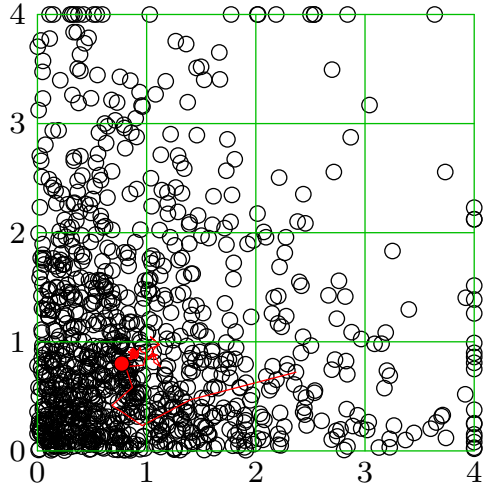}
\end{center}
\caption{Median of a non uniform measure on the unit square in the
plane}
\end{figure}

\subsubsection{Medians in the Poincar\'{e} disc}
\noindent In the two figures below, $M$ is the Poincar\'{e} disc,
the blue points are data points and the red path represents one
trajectory of the inhomogeneous Markov chain $(X_k)_{k\ge 0}$
corresponding to $p=1$, with linear interpolation between the
different steps. The green points are medians computed by the
subgradient method developed in section \ref{rme}.
\begin{figure}[!htpb]
\begin{center}
\includegraphics*[bb= 0 0 376 373, width =2.9in]{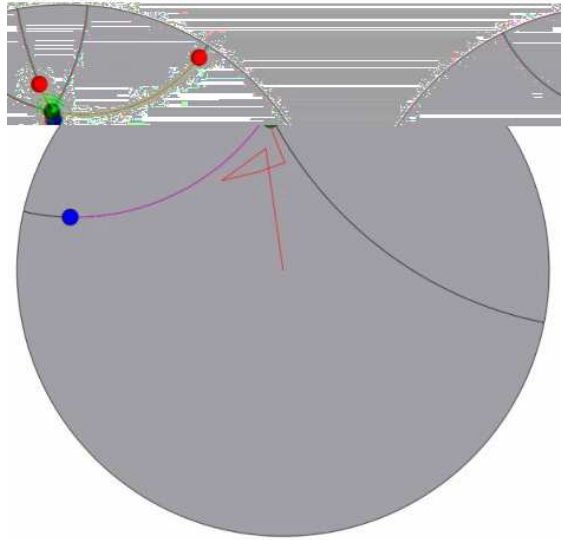}
\end{center}
\caption{Median of three points in the Poicar\'{e}
disc}
\end{figure}
\begin{figure}[!htpb]
\begin{center}
\includegraphics*[bb= 0 0 376 372, width =2.9in]{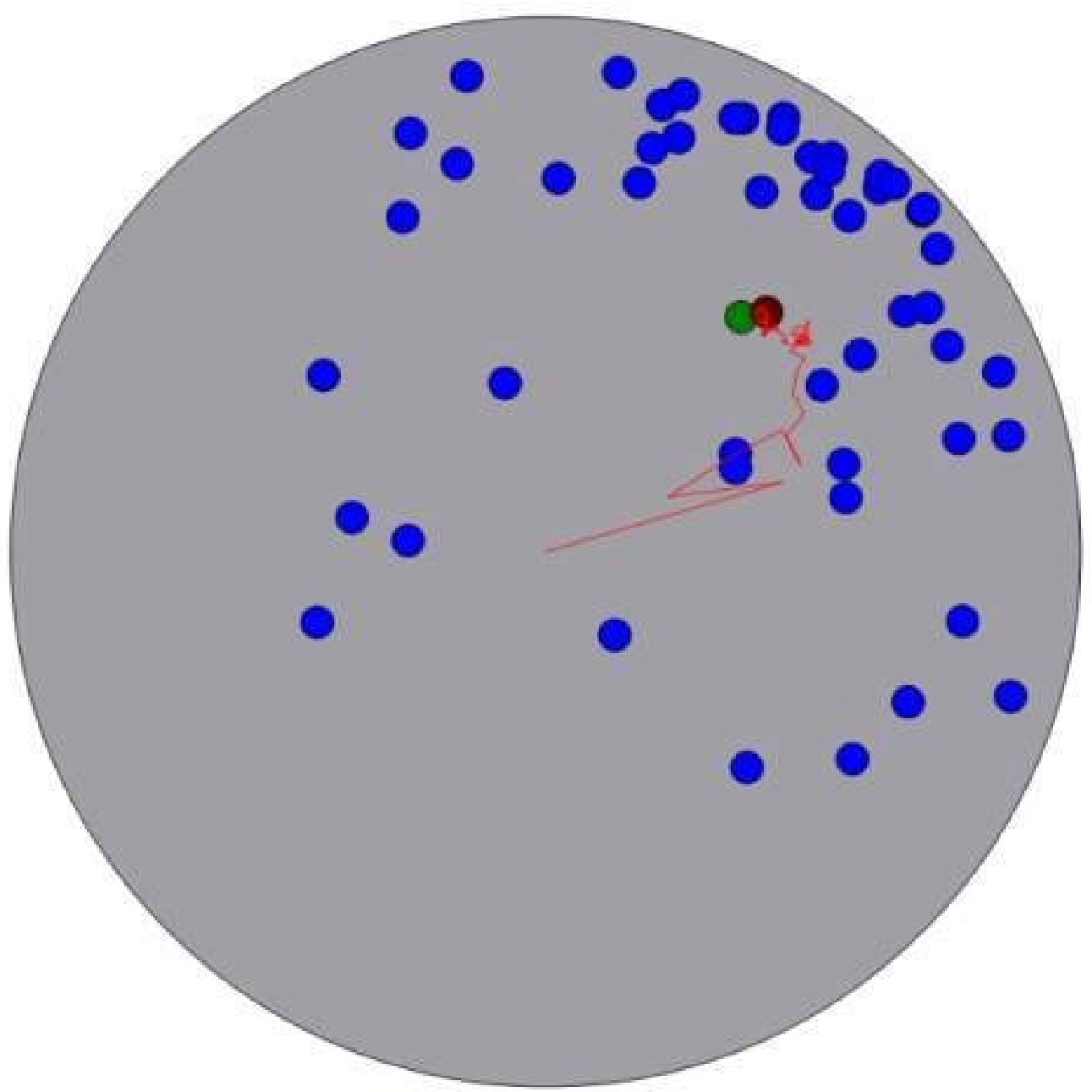}
\end{center}
\caption{Median of points in the Poicar\'{e}
disc}
\end{figure}

\subsection{Computing $p$-means by gradient descent}
Gradient descent algorithms for computing $e_p$ are given in the
following theorem. In view of Theorem \ref{1 convergence}, it
suffices to consider the case when $p>1$.

\begin{theorem}\label{1 gradient descent}
Assume that $p>1$. Let $x_0\in\bar{B}(a,r)$ and for $k\geq0$ define
$$x_{k+1}=\exp_{x_k}(-t_k\grad_{x_k}H_p),$$
where $(t_k)_k$ is a sequence of real numbers such that
$$0<t_k\leq \frac{p\varepsilon^{p+1}}{\pi p^2(2r)^{2p-1}\beta\coth(2\beta r)+p\varepsilon^{p}},\quad\lim_{k\rightarrow\infty}t_k=0\quad\text{and}\quad\sum_{k=0}^{\infty}t_k=+\infty.$$
Then the sequence $(x_k)_k$ is contained in $\bar{B}(a,\rho)$ and
converges to $ e_p$.
\end{theorem}

The following proposition gives the error estimations of the
gradient descent algorithms in Theorem \ref{1 gradient descent}.

\begin{proposition}
Assume that $t_k<C_{p,\mu,K}^{-1}$ for every $k$ in Theorem \ref{1
gradient descent}, then the following
error estimations hold:\\
i) if $1<p<2$, then for $k\geq1$,
\begin{align*}
\rho^2(x_k,e_p)\leq\,&4r^2\prod_{i=0}^{k-1}(1-C_{p,\mu,K}t_i)\\
&+C(\beta,r,p)\bigg(\sum_{j=1}^{k-1}t_{j-1}^2\prod_{i=j}^{k-1}(1-C_{p,\mu,K}t_i)+t_{k-1}^2\bigg):=
b_k;
\end{align*}
ii) if $p\geq2$, then for $k\geq1$,
\begin{align*}
H_p(x_k)-H_p(e_p)\leq\,&(2r)^p\prod_{i=0}^{k-1}(1-C_{p,\mu,K}t_i)\\
&+C(\beta,r,p)\bigg(\sum_{j=1}^{k-1}t_{j-1}^2\prod_{i=j}^{k-1}(1-C_{p,\mu,K}t_i)+t_{k-1}^2\bigg):=
c_k,
\end{align*}
where the constant
$$
C(\beta,r,p)=
\begin{cases}
p^2(2r)^{2p-1}\beta\coth(2\beta r),\,\,& \text{if \,\,$1<p<2$;}\\
p^3(2r)^{3p-4}\left(2\b r\coth (2\b r)+p-2\right),\,\,& \text{if
\,\,$p\geq2$.}
\end{cases}
$$
Moreover, the sequences $(b_k)_k$ and $(c_k)_k$ both tend to zero.
\end{proposition}

\section{Riemannian geometry of Toeplitz covariance matrices\\and
applications to radar target detection}\label{Toeplitz}

In this section we study the Riemannian geometry of the manifold of
Toeplitz covariance matrices of order $n$. The explicit expression
of the reflection coefficients reparametrization and its inverse are
obtained. With the Riemannian metric given by the Hessian of a
K\"{a}hler potential, we show that the manifold is in fact a
Cartan-Hadamard manifold with lower sectional curvature bound $-4$.
The geodesics in this manifold are also computed. Finally, we apply
the subgradient algorithm introduced in section \ref{rme} and the
Riemannian geometry of Toeplitz covariance matrices to radar target
detection. We refer to \cite{Yang4} for more mathematical details of
this section.

\subsection{Reflection coefficients parametrization}
Let $\mathcal{T}_n$ be the set of Toeplitz Hermitian positive
definite matrices of order $n$. It is an open submanifold of
$\mathbf{R}^{2n-1}$. Each element $R_n\in\mathcal{T}_n$ can be
written as
\[
R_n=\left[
\begin{matrix}

r_0                 &     \overline{r}_1      &   \ldots           & \overline{r}_{n-1}\\
r_1      &     r_0      &   \ldots           & \overline{r}_{n-2}\\
\vdots              &     \ddots   &   \ddots           & \vdots\\
r_{n-1}  &     \ldots   &   r_1   & r_0

\end{matrix}
\right]\,.
\]
For every $1\leq k\leq n-1$, the upper left $(k+1)$-by-$(k+1)$
corner of $R_n$ is denoted by $R_k$. It is associated to a $k$-th
order autoregressive model whose Yule-Walker equation is
\[
\left[
\begin{matrix}

r_0                 &     \overline{r}_1      &   \ldots           & \overline{r}_{k}\\
r_1      &     r_0      &   \ldots           & \overline{r}_{k-1}\\
\vdots              &     \ddots   &   \ddots           & \vdots\\
r_{k}  &     \ldots   &   r_1   & r_0

\end{matrix}
\right] \left[
\begin{matrix}
1\\
a_1^{(k)}\\
\vdots\\
a_{k}^{(k)}
\end{matrix}
\right] = \left[
\begin{matrix}
P_k\\
0\\
\vdots\\
0
\end{matrix}
\right]\,,
\]
where $a_1^{(k)},\dots,a_k^{(k)}$ are the optimal prediction
coefficients and $P_{k}=\det R_{k+1}/\det R_k$ is the mean squared
error.

The last optimal prediction coefficient $a_k^{(k)}$ is called the
$k$-th reflection coefficient and is denoted by $\mu_{k}$. It is
easily seen that $\mu_1,\dots,\mu_{n-1}$ are uniquely determined by
the matrix $R_n$. Moreover, the classical Levinson's recursion gives
that $|\mu_k|<1$. Hence, by letting $P_0=r_0$, we obtain a map
between two submanifolds of
$\mathbf{R}^{2n-1}$:\\
\[\varphi:\quad\mathcal{T}_n\longrightarrow \mathbf{R}^*_+\times \mathbf{D}^{n-1},\quad
R_n\longmapsto(P_0,\mu_1,\dots,\mu_{n-1}),\]\\
where $\mathbf{D}=\{z\in\mathbf{C}:|z|<1\}$ is the unit disc of the
complex plane.

Using the Cramer's rule and the method of Schur complement we get
the following proposition.

\begin{proposition}
$\varphi$ is a diffeomorphism, whose explicit expression is
\[\mu_k=(-1)^k \frac{\det S_k}{\det R_k},\quad \text{where}
\quad S_k=R_{k+1}\binom{2,\dots,k+1}{1,\dots,k}\] is the submatrix
of $R_{k+1}$ obtained by deleting the first row and the last column.
On the other hand, if
$(P_0,\mu_1,\dots,\mu_{n-1})\in\mathbf{R}^*_+\times
\mathbf{D}^{n-1}$, then its inverse image $R_n$ under $\varphi$ can
be calculated by the following algorithm:
\[r_0=P_0, \quad r_1=-P_0\mu_1,\]
\[r_k=-\mu_kP_{k-1}+\alpha^T_{k-1}J_{k-1}R_{k-1}^{-1}\alpha_{k-1},\quad 2\leq k \leq n-1,\]
where
\[
\alpha_{k-1}= \left[
\begin{matrix}
r_1\\
\vdots\\
r_{k-1}
\end{matrix}
\right],\quad J_{k-1}= \left[
\begin{matrix}
       0 & \ldots & 0 & 1\\
       0 & \ldots & 1 & 0\\
         & \ldots\\
       1 & \ldots & 0 & 0
\end{matrix}
\right]\quad\text{and}\quad
P_{k-1}=P_0\prod_{i=1}^{k-1}(1-|\mu_i|^2).\]
\end{proposition}

\subsection{Riemannian geometry of Toeplitz covariance matrices}

From now on, we regard $\mathcal{T}_n$ as a Riemannian manifold
whose metric, which is introduced in \cite{Barbaresco2} by the
Hessian of the K\"{a}hler potential $$\Phi(R_n)=-\ln(\det
R_n)-n\ln(\pi e),$$ is given by
\begin{equation}\label{1 Bergman metric}
ds^2=n\frac{dP_0^2}{P_0^2}+\sum_{k=1}^{n-1}(n-k)\frac{|d\mu_k|^2}{(1-|\mu_k|^2)^2},
\end{equation}
where $(P_0,\mu_1,\dots,\mu_{n-1})=\varphi(R_n)$.

The metric \eqref{1 Bergman metric} is a Bergman type metric and it
has be shown in \cite{Yang4} that this metric is not equal to the
Fisher information metric of $\mathcal{T}_n$. But J. Burbea and C.
R. Rao have proved in \cite[Theorem 2]{Burbea} that the Bergman
metric and the Fisher information metric do coincide for some
probability density functions of particular forms. A similar
potential function was used by S. Amari in \cite{Amari} to derive
the Riemannian metric of multi-variate Gaussian distributions by
means of divergence functions. We refer to \cite{Shima} for more
account on the geometry of Hessian structures.

With the metric given by \eqref{1 Bergman metric} the space
$\mathbf{R}^*_+\times \mathbf{D}^{n-1}$ is just the product of the
Riemannian manifolds $(\mathbf{R}^*_+,ds^2_0)$ and
$(\mathbf{D},ds^2_k)_{1\leq k\leq n-1}$, where
$$ds_0^2=n\frac{dP^2_0}{P_0^2}\quad\text{and}\quad
ds^2_k=(n-k)\frac{|d\mu_k|^2}{(1-|\mu_k|^2)^2}.$$ The latter is just
$n-k$ times the classical Poincar\'{e} metric of $\mathbf{D}$. Hence
$(\mathbf{R}^*_+\times \mathbf{D}^{n-1},ds^2)$ is a Cartan-Hadamard
manifold whose sectional curvatures $K$ verify $-4\leq K\leq 0$. The
Riemannian distance between two different points $x$ and $y$ in
$\mathbf{R_+^*}\times \mathbf{D}^{n-1}$ is given by
\[d(x,y)=\bigg(n\sigma(P,Q)^2+\sum_{k=1}^{n-1}(n-k)\tau(\mu_k,\nu_k)^2\bigg)^{1/2},\]
where $x=(P,\mu_1,\ldots,\mu_{n-1})$, $
y=(Q,\nu_1,\ldots,\nu_{n-1})$,
\[\sigma(P,Q)=|\ln(\frac{Q}{P})|\quad\text{and}\quad\tau(\mu_k,\nu_k)=
\frac{1}{2}\ln\frac{1+|\frac{\nu_k-\mu_k}{1-\bar{\mu}_k\nu_k}|}{1-|\frac{\nu_k-\mu_k}{1-\bar{\mu}_k\nu_k}|}.\]

The geodesic from $x$ to $y$ in $\mathcal{T}_n$ parameterized by arc
length is given by
\[\gamma(s,x,y)=(\gamma_0(\frac{\sigma(P,Q)}{d(x,y)}s),\gamma_1(\frac{\tau(\mu_1,\nu_1)}{d(x,y)}s),\ldots,\gamma_1(\frac{\tau(\mu_{n-1},\nu_{n-1})}{d(x,y)}s)),\]
where $\gamma_0$ is the geodesic in $(\mathbf{R_+^*},ds^2_0)$ from
$P$ to $Q$ parameterized by arc length and for $1\leq k\leq n-1$,
$\gamma_k$ is the geodesic in $(\mathbf{D},ds^2_k)$ from $\mu_k$ to
$\nu_k$ parameterized by arc length. More precisely,
\[\gamma_0(t)=Pe^{t\sign(Q-P)},\]
and for $1\leq k\leq n-1$,
\[\gamma_k(t)=
\frac{(\mu_k+e^{i\theta_k})e^{2t}+(\mu_k-e^{i\theta_k})}{(1+\bar{\mu}_ke^{i\theta_k})e^{2t}+(1-\bar{\mu}_ke^{i\theta_k})}
,\quad\text{where}\quad\theta_k=\arg\frac{\nu_k-\mu_k}{1-\bar{\mu}_k\nu_k}.\]
Particularly,
\[\gamma^{\prime}(0,x,y)=(\gamma_0^{\prime}(0)\frac{\sigma(P,Q)}{d(x,y)},\gamma_1^{\prime}(0)\frac{\tau(\mu_1,\nu_1)}{d(x,y)},\ldots,
\gamma_{n-1}^{\prime}(0)\frac{\tau(\mu_{n-1},\nu_{n-1})}{d(x,y)}).\]

Let $v=(v_0,v_1,\ldots,v_{n-1})$ be a tangent vector in
$T_x(\mathbf{R}_+^*\times \mathbf{D}^{n-1})$, then the geodesic
starting from $x$ with velocity $v$ is given by
\[\zeta(t,x,v)=(\zeta_0(t),\zeta_1(t),\ldots,\zeta_{n-1}(t)),\]
where $\zeta_0$ is the geodesic in $(\mathbf{R}_+^*,ds^2_0)$
starting from $P$ with velocity $v_0$ and for $1\leq k\leq n-1$,
$\zeta_k$ is the geodesic in $(\mathbf{D},ds^2_k)$ starting from
$\mu_k$ with velocity $v_k$. More precisely,
\[\zeta_0(t)=Pe^{\frac{v_0}{P}t},\]
and for $1\leq k\leq n-1$,
\[\zeta_k(t)=\frac{(\mu_k+e^{i\theta_k})e^{\frac{2|v_k|t}{1-|\mu_k|^2}}+(\mu_k-e^{i\theta_k})}{(1+\bar{\mu}_ke^{i\theta_k})e^{\frac{2|v_k|t}{1-|\mu_k|^2}}+(1-\bar{\mu}_ke^{i\theta_k})}
,\quad\text{where}\quad\theta_k=\arg v_k.\]

\subsection{Radar simulations}
Now we give some simulating examples of the median method applied to
radar target detection.

Since the autoregressive spectra are closely related to the speed of
targets, we shall first investigate the spectral performance of the
median method. In order to illustrate the basic idea, we only
consider the detection of one fixed direction. The range along this
direction is subdivided into 200 lattices in which we add two
targets, the echo of each lattice is modeled by an autoregressive
process. The following Figure \ref{initial spectra} gives the
initial spectra of the simulation, where $x$ axis represents the
lattices and $y$ axis represents frequencies. Every lattice is
identified with a $1\times8$ vector of reflection coefficients which
is calculated by using the regularized Burg algorithm
\cite{Barbaresco5} to the original simulating data. The spectra are
represented by different colors whose corresponding values are
indicated in the colorimetric on the right.

\begin{figure}[!htpb]
\begin{center}
\includegraphics*[bb= -181   102   778   740, width =3.8in]{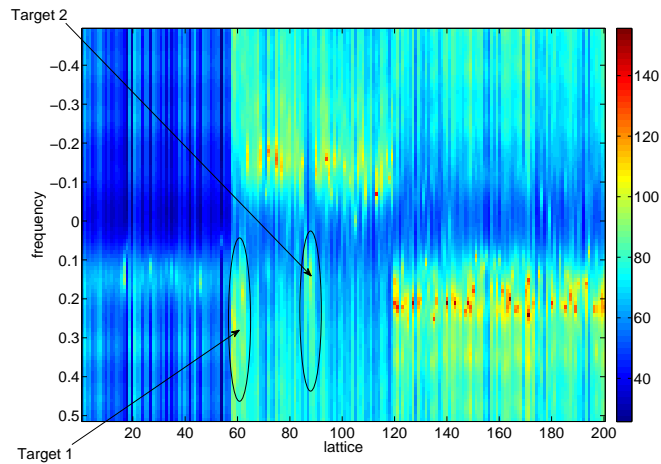}
\end{center}
\caption{Initial spectra with two added targets}\label{initial
spectra}
\end{figure}

\begin{figure}[!htpb]
\begin{center}
\includegraphics*[bb= -211    71   806   769, width =3.75in]{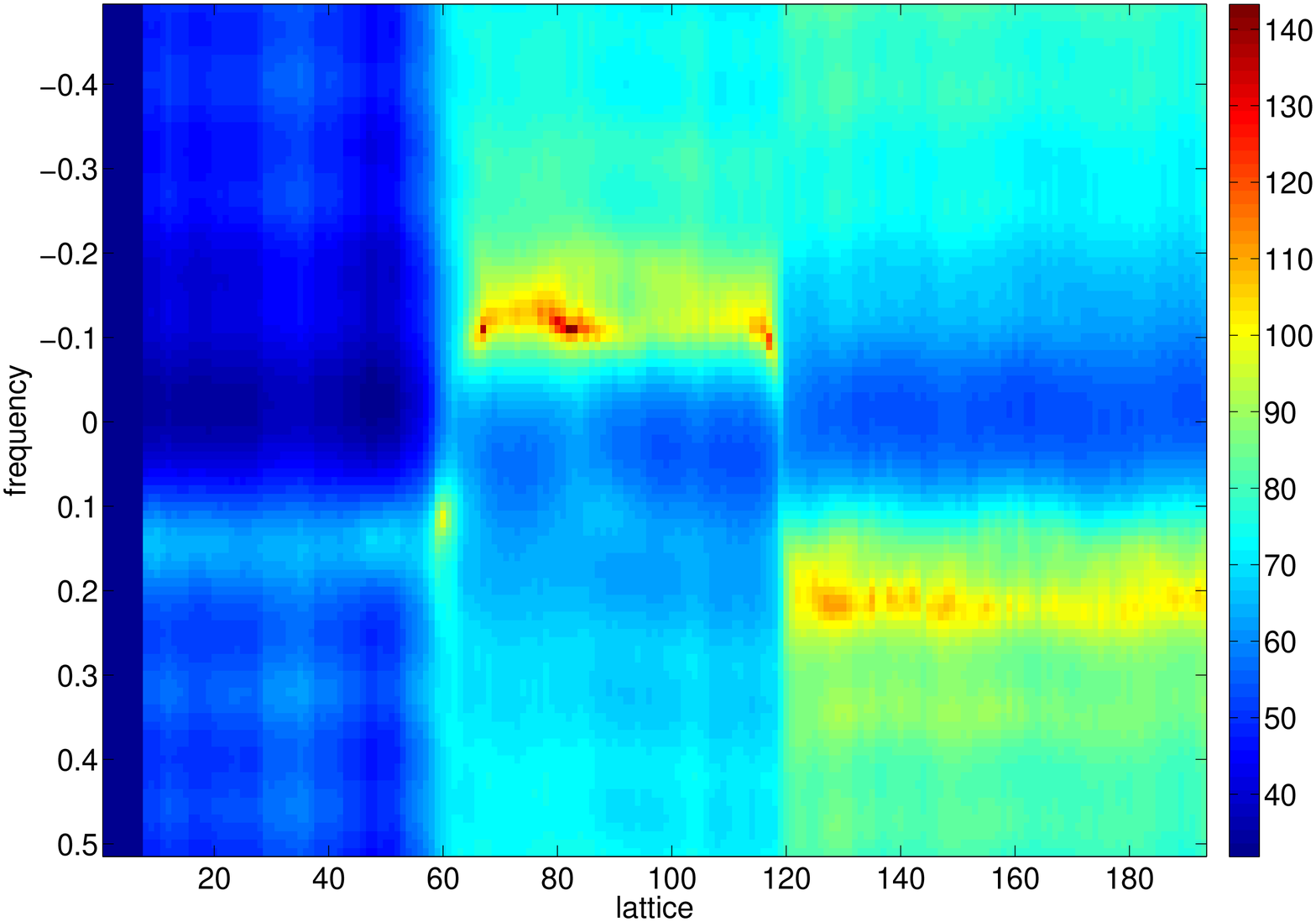}
\end{center}
\caption{Median spectra}\label{Median spectra}
\begin{center}
\includegraphics*[bb= -181   102   778   740, width =3.75in]{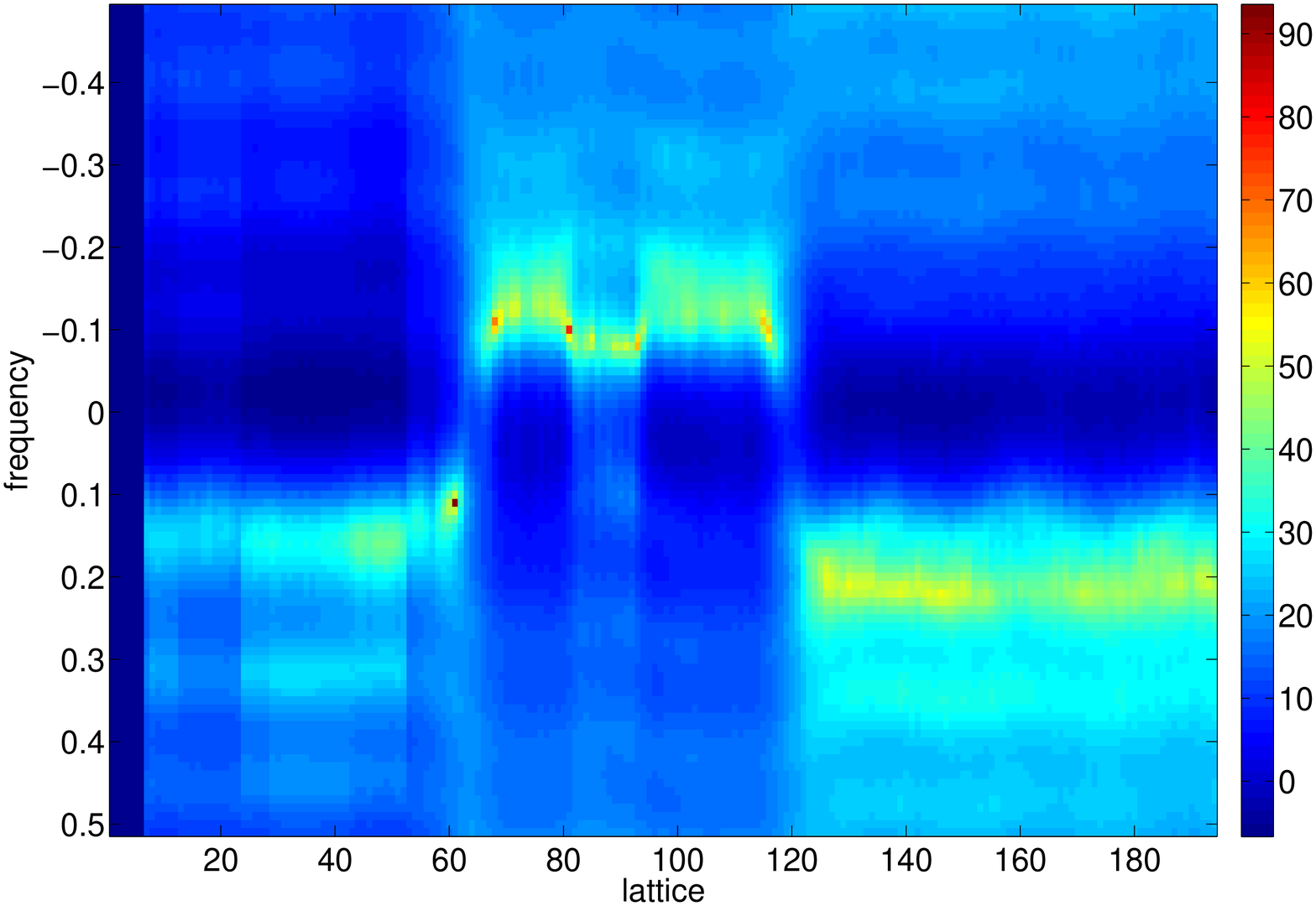}
\end{center}
\caption{Barycenter spectra}\label{Barycenter spectra}
\end{figure}

\newpage
For every lattice, by using the subgradient algorithm, we calculate
the median of the window centered on it and consisting of 15
lattices and then we get the spectra of medians shown in Figure
\ref{Median spectra}. Furthermore, by comparing it with Figure
\ref{Barycenter spectra} which are spectra of barycenters, we see
that in the middle of the barycenter spectra, this is just the place
where the second target appears, there is an obvious distortion.
This explains that median is much more robust than barycenter when
outliers come.

The principle of target detection is that a target appears in a
lattice if the distance between this lattice and the median of the
window around it is much bigger than that of the ambient lattices.
The following Figure \ref{Detection by median} shows that the two
added targets are well detected by the median method, where $x$ axis
represents lattice and $y$ axis represents the distance in
$\mathcal{T}_8$ between each lattice and the median of the window
around it.

\begin{figure}[!htpb]
\begin{center}
\includegraphics*[bb=  -181   102   778   740, width =3.6in]{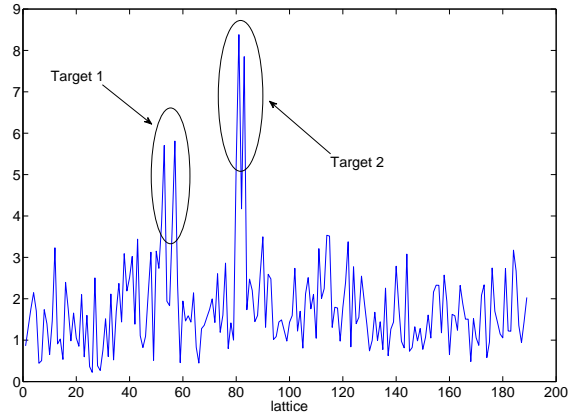}
\end{center}
\caption{Detection by median}\label{Detection by median}
\end{figure}

We conclude our discussion by showing the performance of the median
method in real target detection. As above, we give the images of
autoregressive spectra and the figure of target detection obtained
by using real data which are records of a radar located on a coast.
These records consist of about 5000 lattices of a range of about
10km-140km as well as 109 azimuth values corresponding to
approximately 30 scanning degrees of the radar. For simplicity we
consider the data of all the lattices but in a fixed direction,
hence each lattice corresponds to a $1\times8$ vector of reflection
coefficients computed by applying the regularized Burg algorithm to
the original real data. Figure \ref{Initial spectra of real radar
data} gives the initial autoregressive spectra whose values are
represented by different color according to the colorimetric on the
right. For each lattice, by using the subgradient algorithm, we
calculate the median of the window centered on it and consisting of
17 lattices and then we get the spectra of medians shown in Figure
\ref{Median spectra of real radar data}.

\begin{figure}[!htpb]
\begin{center}
\includegraphics*[bb=  -181   102   778   740, width =3.8in]{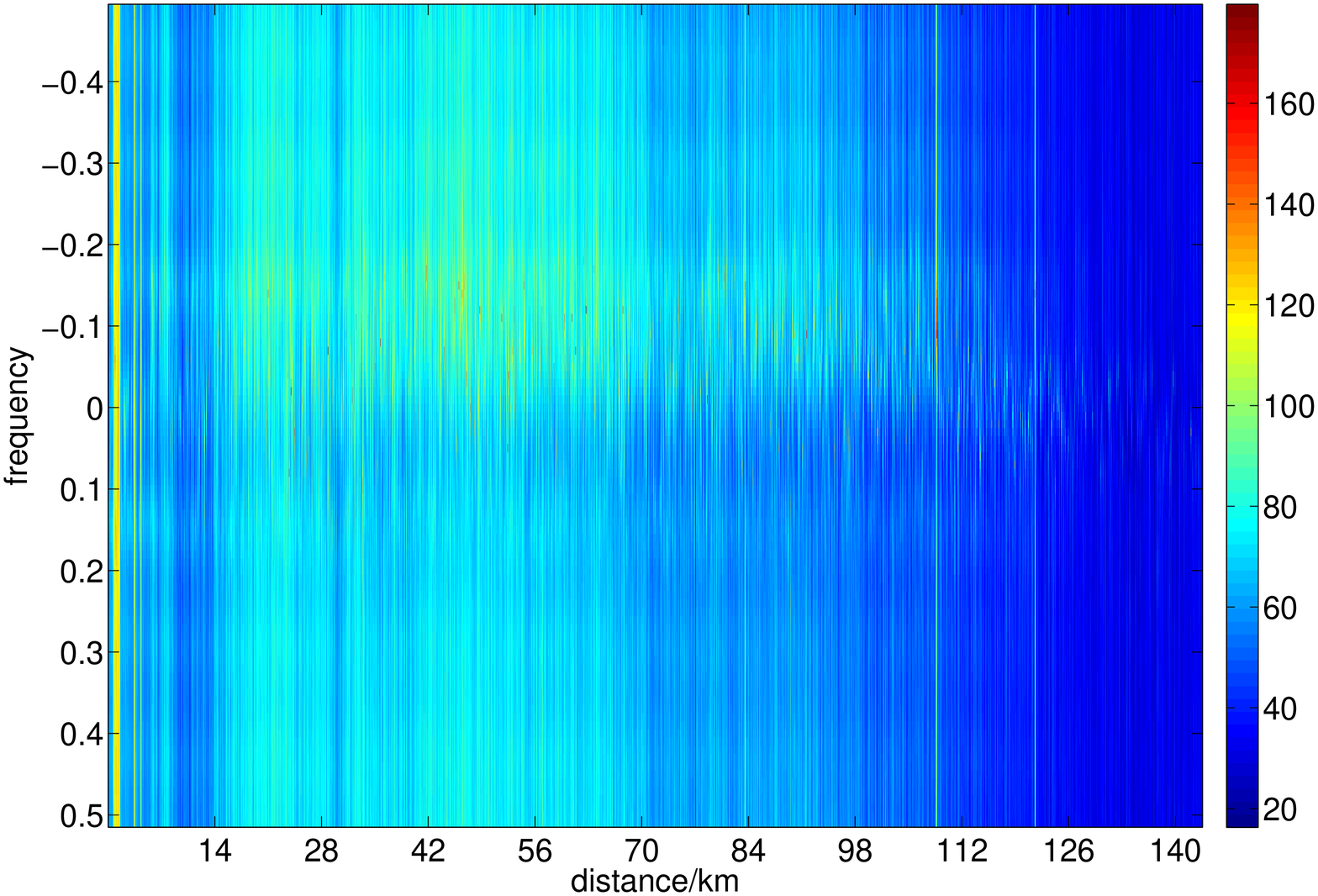}
\end{center}
\caption{Initial spectra of real radar data}\label{Initial spectra
of real radar data}
\begin{center}
\includegraphics*[bb=  -181   102   778   740, width =3.8in]{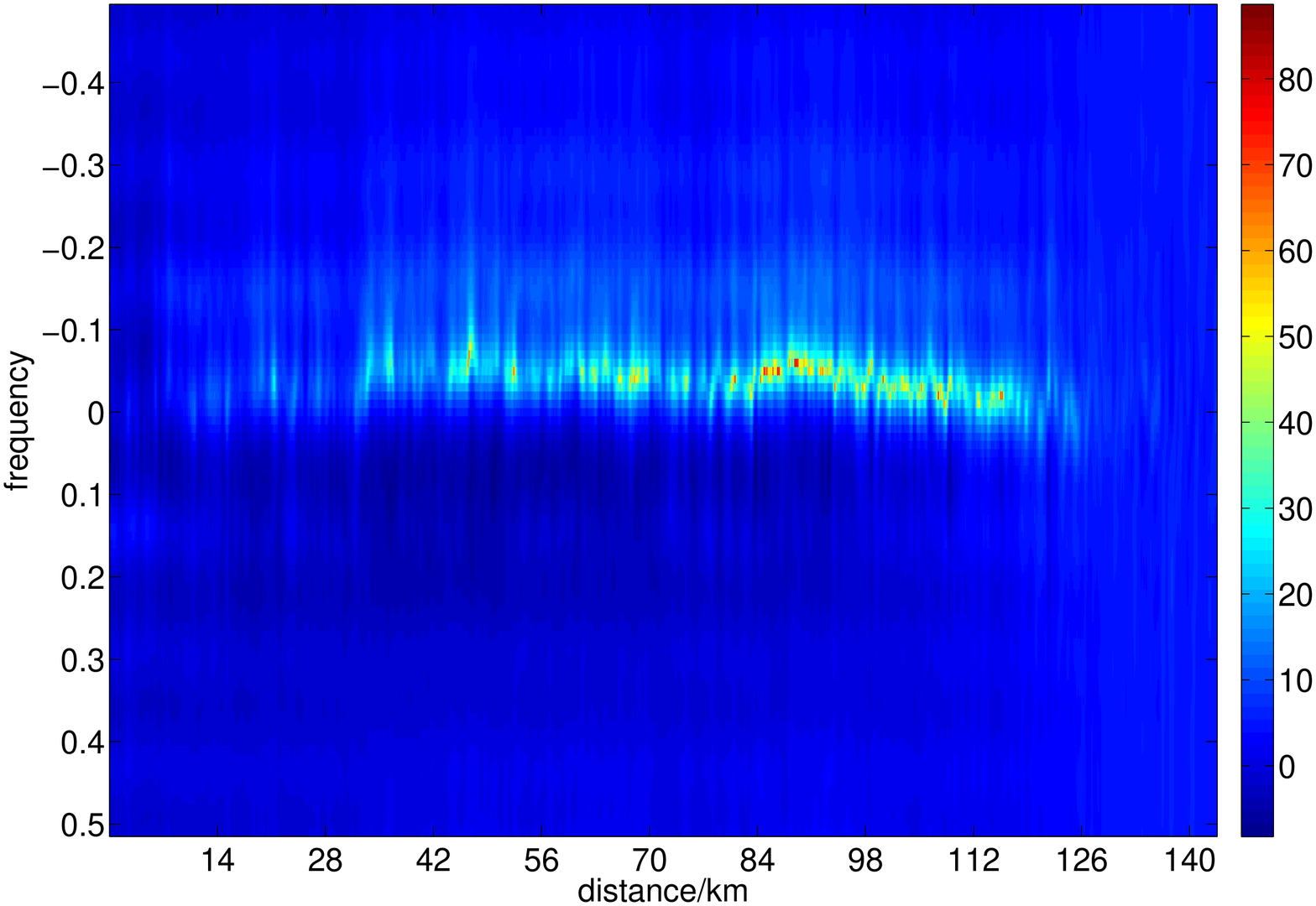}
\end{center}
\caption{Median spectra of real radar data}\label{Median spectra of
real radar data}
\end{figure}

\begin{figure}[!htpb]
\begin{center}
\includegraphics*[bb= -211    71   806   769, width =3.8in]{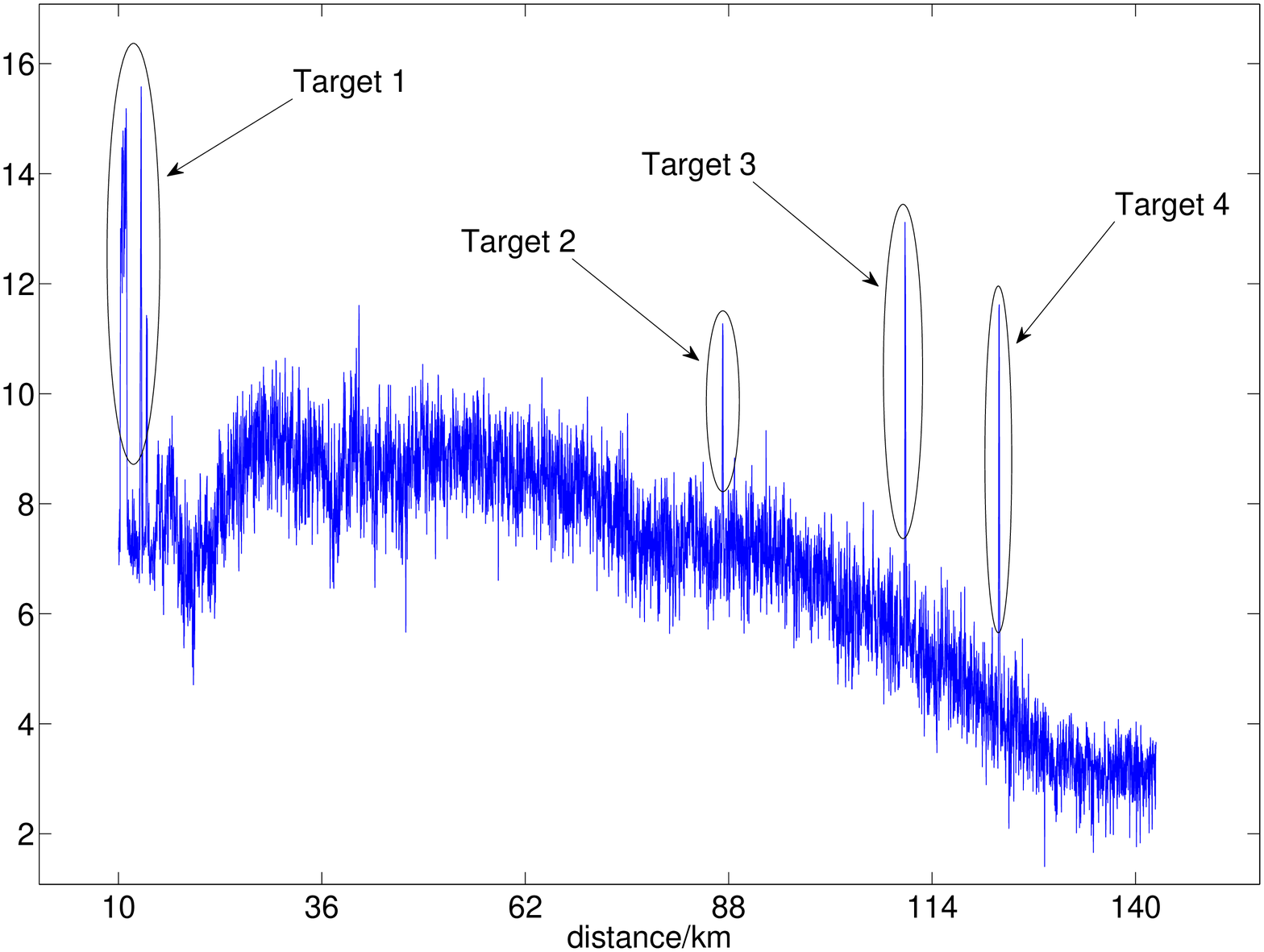}
\end{center}
\caption{Real detection by median}\label{Real detection by median}
\end{figure}

In order to know in which lattice target appears, we compare the
distance between each lattice and the median of the window around
it. The following Figure \ref{Real detection by median} shows that
the four targets are well detected by our method, where $x$ axis
represents distance and $y$ axis represents the distance in
$\mathcal{T}_8$ between each lattice and the median of the window
around it.


\end{document}